  \documentclass[a4paper,11pt]{article}
    \usepackage[dvips]{graphicx}
  \usepackage[]{fancyheadings}
  \usepackage{makeidx,epsf}
  \usepackage{latexsym}
  \usepackage{hyperref,nameref} 
  \usepackage{url}
  \scrollmode
  \makeindex
\font\tenbit=cmmib10 at 12pt

\def\bzxplain{{\hbox{\tenbit \char"1F}}}%   chi
\def\bzx{{\raise2pt\bzxplain}}
\def\normalaccent{!}
\def\factorial {!}
\catcode`\!=\active
\def\boldmath#1{\relax\ifmmode{\bf #1}\else\normalaccent\fi}
\def\!{\ifmmode\mskip-\thinmuskip\else\char'74\fi}

\def\zx{\raise2pt\hbox{$\chi$}}

\def\prt{\partial}
\def\dd#1#2{{\partial #1\over \partial #2}}
\def\dep#1.#2.{\dd{#1}{#2}}
\def\depa#1.#2.#3.{\left.\dd{#1}{#2}\right|_{#3}}
\def\dtot#1.#2.{{d #1\over d #2}}
\def\deps#1.#2.{\left(\prt #1/\prt #2\right)}
\def\dtots#1.#2.{\left(d #1/d #2\right)}
\def\at#1\at#2.{\left.#1\right|_{#2}}

\def\di#1.{\,\,d#1}
\def\D#1.{\Delta\!#1}

\def\AeAe1{A^+ ~ e^{imx} + A^- \, e^{-imx}}
\def\AeAe2{A^+ ~ e^{imx} - A^- \, e^{-imx}}
\def\AeAe#1{A^+ \, e^{imx} #1 A^- \, e^{-imx}}

\def\Dph1{\phi'}
\def\Det1{\eta'}

\def\rpi4{{1 \over 4 \pi}}

\def\SO3{{\left [ \matrix {  C & S  \cr -S & C \cr } \right ] }}

\let!=\boldmath

  \advance\voffset by  -1.40 truecm      % (-4.80 two / -3.40 one )
  \advance\hoffset by  -1.40 truecm      % ( 1.50 two / -1.00 one )
\newcommand{\ltitle}[1]{{\normalsize {\bf #1 }}} 

\def\thehead$#1${#1}

\begin{document}

\title { {\rm \normalsize
A semi-numerical computation for the added mass coefficients of \\ %[+0.025cm]
an oscillating hemi-sphere at very low and very high frequencies }
       }
\author{ \normalsize
     $\hbox{M.A. Storti}^1$ 
     and 
     $\hbox{J. D'El\'{\i}a}^1$
         \\ [-0.02 cm] \normalsize
     1: Centro Internacional de M\'etodos 
     Computacionales en Ingenier\'{\i}a (CIMEC)
         \\ [-0.02 cm] \normalsize
     INTEC (CONICET-UNL), G\"uemes 3450, 3000-Santa Fe, Argentina \\
  \url{mstorti@intec.unl.edu.ar},
  \url{jdelia@intec.unl.edu.ar} \\
  \url{http://venus.arcride.edu.ar/CIMEC/}
       }
\date{\empty}
\maketitle
% -----------------------------------------------------------------

\pagenumbering{arabic} 
\normalsize 
\begin {abstract}  
{\baselineskip 18 pt % interline 14 - 20
 \noindent
 A floating hemisphere under forced harmonic oscillation at very 
 high and very low frequencies is considered. The problem is reduced 
 to an elliptic one, that is, the Laplace operator in the exterior 
 domain with standard Dirichlet and Neumann boundary conditions, so 
 the flow problem is simplified to standard ones, with well known 
 analytic solutions in some cases. The general procedure is based 
 in the use of spherical harmonics and its derivation is based on 
 a physics insight. The results can be used to test the accuracy 
 achieved by numerical codes as, for example, by finite elements 
 or boundary elements.
\par} 
\end{abstract}
\tableofcontents
\section{Classification (AMS/MSC)}

\begin{itemize}
\item[\tt 76Bxx] Incompressible inviscid fluids
\item[\tt 76B07] Free-surface potential flows
\item[\tt 33C55] Spherical harmonics
\item[\tt 33F05] Numerical approximation
\end{itemize}

\section{Introduction}

A semi-numerical computation of the added mass for a floating 
hemisphere is shown, under an harmonic forced oscillation on 
the free surface of an irrotational, incompressible fluid, and
linearized boundary conditions. Two standard problems are
considered: the heave and the surge modes, that is, vertical 
and horizontal oscillations, respectively. For simplicity, our 
attention is restricted to a fluid of infinite depth. The added 
mass coefficients at very-high and very-low frequencies found 
in this work are in agreement with those given in literature 
by another strategies.

\noindent
As it is noted by Hulme \cite{rf:hulme}, the hydrodynamic 
formulation of a floating hemisphere is analogous to the 
two-dimensional circular cylinder. The added mass coefficients 
are computed as $A_{kk}'=A_{kk}/(\rho V)$, and the damping ones, 
as $D_{kk}'=D_{kk}/(\rho V \omega)$, where $V=(2/3) \pi R^3$ is 
the hemisphere volume, $\rho$ is fluid density, and $\omega$ is 
the circular frequency of the oscillation. The asymptotic values 
of these coefficients, for very slow and very high frequencies, 
can be obtained by analytical calculus, for instance, by a 
variable separation or image methods. 
For the surge/sway mode at very slow frequency, the boundary 
condition $\phi_{,z}=0$, where $\phi$ is the velocity potential, 
is equivalent to a symmetry operation respect the plane $z=0$ 
and, then, corresponds to the solution of a sphere oscillating 
in an infinity medium. The added mass for the last case is half 
of the displaced volume, e.g. see \cite{rf:landweber}, then, 
the surge/sway added mass coefficient is $A'_{11}= 1/2$, 
respect to the true displaced mass $ (2/3)\pi R^3\rho$, where 
the half factor is due to the analytic prolongation. On the 
other hand, the asymptotic values of the added mass in heave 
mode are not too easy to obtain, and they may be obtained by 
a semi-numerical computation with spherical harmonics.

\noindent
The reason for doing this work is twofold. First, the method of 
solution adopted here is based by a physics insight rather a 
mathematical one. Second, the computation is near-exact, so the 
results can be used to test the accuracy achieved by numerical 
codes, as finite-element and boundary-element ones, adapted to 
wave-drag and seakeeping flow problems, e.g. see 
\cite{rf:asmejfe, rf:inmf, rf:dnlabso, rf:dnlmar2, rf:ohkusu, 
rf:jdeliathesis}.

\section{An oscillating hemisphere}

An oscillating hemisphere in a forced motion is considered. 
The flat face of the hemisphere is on the free surface of an
irrotational and incompressible fluid, without a mean flow, 
where the fluid depth is assumed as infinity. The upward 
direction is $z$ and the hydrostatic equilibrium plane is 
$z=0$. Due to the symmetry, a spherical coordinate system 
is chosen as
\begin{equation}
   \cases {
    z = r \cos \theta           & ; \cr
    x = r \sin \theta \cos \phi & ; \cr
    y = r \sin \theta \sin \phi & .
          }
\label{eq-nobl150}
\end{equation}
At very high frequencies, the free surface boundary condition
shrinks to the simple form $\phi=0$. But, by anti-symmetry 
respects to the plane $z=0$, the flow problem is equivalent 
to solve the modified one
\begin{equation}
   \cases {
   \Delta \phi = 0            & in $ \Omega'   $ ; \cr
   \phi_{,n}   = \cos \theta  & at $ \Gamma'_e $ ; \cr
          \phi = 0            & at $ \Gamma'_0 $ ;
          } 
\label{eq-nobl170}
\end{equation}
where $\phi$ is the velocity potential, $\Delta$ is the Laplace 
operator, $\Omega'$, $\Gamma'_e$ and $\Gamma'_0$ are the extended 
flow domain, extended hemisphere surface and extended free surface, 
respectively, through the reflection plane $z=0$. As the free 
surface boundary condition for very high frequencies is $\phi=0$, 
then, its right hand side term has been extended in an anti-symmetric 
way. That is, the extended problem at very high frequencies is the 
same that a sphere in infinite medium. 
On the other hand, the free surface boundary condition for very 
slow frequencies is $\phi_{,n}=0$, and for the extended 
problem, its right hand side must be extended in a symmetrical 
way and, then, 
\begin{equation}
  \cases {
  \Delta \phi = 0               & in $ \Omega'   $ ; \cr
   \phi_{,n}  = | \cos \theta | & at $ \Gamma'_e $ ; \cr
         \phi = 0               & at $ \Gamma'_0 $ ;
          } 
\label{eq-nobl180}
\end{equation}
where, due the module on $|\cos \theta |$, the slow frequency
radiation problem does not have, in general, a closed solution 
and, then, it must be found with another resources like spherical 
harmonics, as it is considered in this work.

\section{The very-high and very-low frequencies limits}

The free-surface boundary condition in the limits of very-low 
and very frequencies
\begin{equation} 
        \phi_n = {\omega^2 \over g} \phi \ ;
\label{eq-add010}
\end{equation}
is reduced to the homogeneous Neumann and Dirichlet boundary 
conditions, respectively, where $g$ is the gravity acceleration.
Also, the radiation boundary condition at infinity imposes that 
the velocity potential $\phi$ tends to zero, so the corresponding 
Partial Differential Equation (PDE) system is
\begin{equation} 
    \cases {
 \Delta \phi     = 0          & in $\Omega  $ ; \cr
        \phi_{,n}= i \omega h & at $\Gamma_B$ ; \cr
        \phi_{,n}= 0          & at $\Gamma_F$ for low frequencies; \cr
        \phi     = 0          & at $\Gamma_F$ for high frequencies; \cr 
        | \phi | \to 0        & for $ | !x | \to \infty $ ;
           }
\label{eq-add020}
\end{equation}
where the load $h$ is the normal displacement of the mode under 
consideration. It can be seen that, under these conditions, the 
flow problem is transformed in a standard elliptic one, whose 
solution is real valued (in reality it is an imaginary one but 
it can be transformed by means of a simple re-definition). It is 
assumed that the load $h$ is real (that is, the body motion is in 
phase). The added mass for the mode motion is found from
\begin{equation} 
  a_{jj} = - {1 \over \omega^2} 
             \int_{\Gamma} d \Gamma i \omega \phi_{,n}
         = - \int_{\Gamma} d \Gamma          \psi_{,n} \ ;
\label{eq-add030}
\end{equation}
where now $\psi = -(i/\omega) \phi $, so in the limits 
$\omega \to 0$ and $\omega \to \infty$, the function $\psi$ is real.

%%%%%%%%%%%%%%%%%%%%%%%%%%%%%%%%%%%%%%%%%%%%%%%%%%%%%%%%%%%%%%%%%%%
\begin{figure}[tb]
\centerline{\includegraphics[width=8.0truecm]{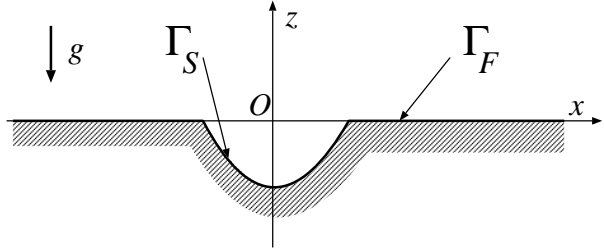}}
\caption{Geometrical description of a seakeeping-like flow problem.}
\label{fg-descrip}
\end{figure}
%%%%%%%%%%%%%%%%%%%%%%%%%%%%%%%%%%%%%%%%%%%%%%%%%%%%%%%%%%%%%%%%%%%

\noindent
By symmetry, the Eqns (\ref{eq-add020} c-d) can be reproduced
extending the flow problem to $z>0$, by means of a mirror
body image and extending the load $h$ in an appropriate way.
For instance, the homogeneous Neumann boundary condition is
obtained extending the load in a symmetrical way with 
respect to $z=0$, that is, 
\begin{equation} 
   h (x,y,z) = + h (z,y,-z) \ ;
\label{eq-add040}
\end{equation}
while the Dirichlet boundary condition can be obtained 
extending in a skew-symmetrical way
\begin{equation} 
   h (x,y,z) = + h (z,y,-z) \ ;
\label{eq-add050}
\end{equation}
For example, the boundary boundary condition for the hemisphere 
in the heave-mode is
\begin{equation} 
   \psi_{,n} = \cos \varphi \ ;
\label{eq-add060}
\end{equation}
where $(r,\varphi,\theta)$ is the spherical coordinate system
with origin in the center of the hemisphere such as
\begin{equation}
    \cases {
    z = r \cos  \varphi                  & ; \cr
    x = r \sin  \varphi  \cos  \theta    & ; \cr
    y = r \sin  \varphi  \sin  \theta    & ;
           }
\label{eq-add070}
\end{equation}
here the Hildebrand's convention \cite{rf:hildebrand} is used. 
Then, the limits at very-low and very-high frequencies for the 
hemisphere can be computed from the sphere ones, but with the 
loads
\begin{equation} 
   \cases {
      \psi_{,n} = | \cos  \varphi |  & for low frequencies  ; \cr
      \psi_{,n} =   \cos  \varphi    & for high frequencies ;
          }
\label{eq-add080}
\end{equation}
%
%%%%%%%%%%%%%%%%%%%%%%%%%%%%%%%%%%%%%%%%%%%%%%%%%%%%%%%%%%%%%%%%%%%
\begin{figure}[tb]
\centerline{\includegraphics[width=13.0truecm]{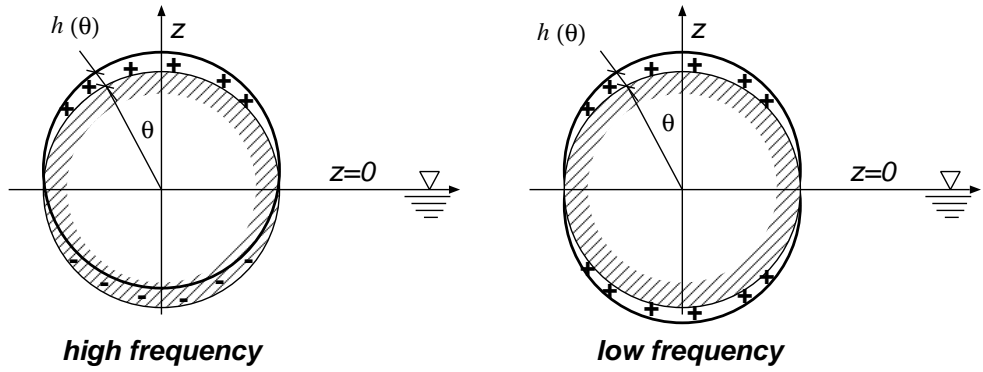}}
\caption{Extension of the load $h$ for the heave-mode at 
         very-low and very-high frequencies.}
\label{fg-heavemo}
\end{figure}
%%%%%%%%%%%%%%%%%%%%%%%%%%%%%%%%%%%%%%%%%%%%%%%%%%%%%%%%%%%%%%%%%%%
%
%%%%%%%%%%%%%%%%%%%%%%%%%%%%%%%%%%%%%%%%%%%%%%%%%%%%%%%%%%%%%%%%%%%
\begin{figure}[tb]
\centerline{\includegraphics[width=6.5truecm]{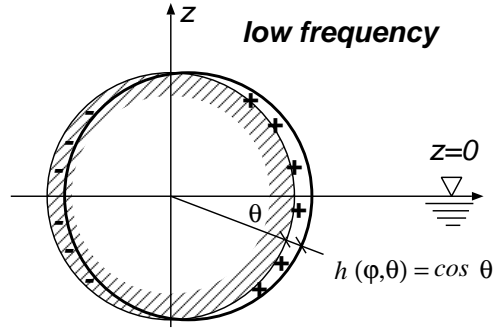}}
\caption{Extension of the load $h=\cos\theta$ 
         for the surge-mode at very-low frequencies.}
\label{fg-surgelf}
\end{figure}
%%%%%%%%%%%%%%%%%%%%%%%%%%%%%%%%%%%%%%%%%%%%%%%%%%%%%%%%%%%%%%%%%%%
%
%%%%%%%%%%%%%%%%%%%%%%%%%%%%%%%%%%%%%%%%%%%%%%%%%%%%%%%%%%%%%%%%%%%
\begin{figure}[tb]
\centerline{\includegraphics[width=9.0truecm]{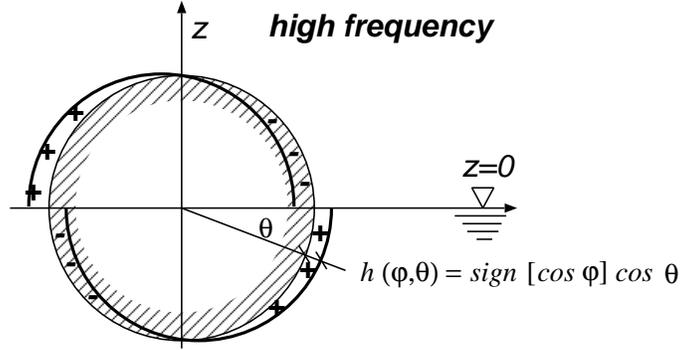}}
\caption{Extension of the load $h=\hbox{sign}~(\cos\varphi)\cos\theta$
         for surge-mode at very-high frequencies.}
\label{fg-surgehf}
\end{figure}
%%%%%%%%%%%%%%%%%%%%%%%%%%%%%%%%%%%%%%%%%%%%%%%%%%%%%%%%%%%%%%%%%%%
%
see figure~\ref{fg-heavemo}. Similarly, for the surge mode 
(oscillation along the $x$-axis) the equivalent load is
\begin{equation} 
   h ~ (\varphi \theta) = \cos \theta \ ;
\label{eq-add090}
\end{equation}
and the extensions for very-low and very-high frequencies are
\begin{equation} 
   \cases {
      \psi_{,n} = \cos  \theta  & for low frequencies ; \cr
      \psi_{,n} = \hbox{sign} ~ ( \cos \varphi ) ~
                  \cos  \theta  & for high frequencies ;
          }
\label{eq-add100}
\end{equation}
\section{Solution of the flow problems}
The previous flow problems can be solved, for instance, by an 
analytical way or by series. The heave solution at very-high 
frequencies and the surge one at very-low frequencies are the 
same of a sphere in an infinity medium and uniform velocity, so 
the additional mass is a half of the displaced mass, that is,
\begin{equation} 
   \cases {
   a_{33} (\omega \to \infty) = \frac \pi 3 \rho R^3 & ; \cr
   a_{11} (\omega \to      0) = \frac \pi 3 \rho R^3 & .
          }   
\label{eq-add110}
\end{equation}
In the other two cases it is necessary to make an expansion of
the sources by means of spherical harmonics.

\section{Spherical harmonics}

The solution of the exterior potential problem
\begin{equation} 
   \cases {
 \Delta \psi = 0                  & for $r > 1$ ; \cr
        \psi = f (\varphi,\theta) & at  $r = 1$ ;
          }
\label{eq-add120}
\end{equation}
where $\psi = \psi (\varphi,\theta)$, it can be solved 
expanding the function $ f (\varphi,\theta)$ in terms 
of the harmonics
\begin{eqnarray} 
  f (\varphi,\theta) = 
  \sum_{n = 0}^{\infty} a_{n0} ~ P_n   (\cos \varphi) +
  \sum_{n = 0}^{\infty} ~ \sum_{m = 1}^n 
  \left [   a_{nm} ~ \cos (m \theta) 
          + b_{nm} ~ \sin (m \theta) \right ] 
                               ~ P_n^m (\cos \varphi) \ ;
\label{eq-add130}
\end{eqnarray}
where
\begin{eqnarray}
  a_{n0} = \frac {2n + 1} {4 \pi} ~ 
           \int_{r = 1} d \Gamma 
           f (\varphi,\theta) P_n (\cos \varphi) \ ;
\label{eq-add140}
\end{eqnarray}
\begin{eqnarray}
  a_{nm} = \frac {2n + 1} {2 \pi} ~
           \frac {(n-m) \factorial} {(n+m)\factorial} ~
           \int_{r = 1} d \Gamma 
           f (\varphi,\theta) 
           P_n^m (\cos \varphi) \cos m \theta  \ ;
\label{eq-add150}
\end{eqnarray}
\begin{eqnarray}
  b_{nm} = \frac{2n+1}{2\pi} ~
           \frac{(n-m) \factorial} {(n+m) \factorial} ~
           \int_{r=1} d \Gamma 
           f (\varphi,\theta) 
           P_n^m (\cos \varphi) \sin m \theta  \ ;
\label{eq-add160}
\end{eqnarray}
where $ d \Gamma = \sin \varphi ~d \theta ~ d \varphi $ is the
solid angle differential in spherical coordinates. Once this 
expansion is computed, the exterior potential can be written as
\begin{eqnarray} 
  \psi ~ (r,\varphi,\theta) = 
  \sum_{n=0}^{\infty} Y_n (\varphi,\theta) ~ r^{-(n+1)} \ ; 
\label{eq-add170}
\end{eqnarray}
where
\begin{equation} 
 Y_n (\varphi,\theta) = 
 a_{n0} P_n (\cos \varphi) + 
 \sum_{m=1}^n \left  [ 
 a_{nm} \cos (m \theta) + b_{nm} \sin (m \theta) 
              \right ]    P_n^m (\cos \varphi) \ ; 
\label{eq-add180}
\end{equation}
now, the Neumann problem can be solved taking derivatives with
respect to $r$ and evaluating at $r=1$, obtaining the expression  
\begin{equation} 
   h (\varphi, \theta) = - 
   \sum_{n=0}^{\infty} (n+1) ~ Y_n (\varphi,\theta) \ ;
\label{eq-add190}
\end{equation}
from which analogous relations are obtained
\begin{eqnarray}
  a_{n0} = \frac {2n+1} {4 \pi (n+1)} \int_{r=1} d \Gamma ~ 
           h (\varphi,\theta) ~ P_n (\cos \varphi)    
           \qquad \qquad ; \\
  a_{nm} = \frac {2n+1} {2 \pi (n+1)} 
           \frac {(n-m) \factorial } {(n + m) \factorial} 
           \int_{r=1} d \Gamma ~
           h (\varphi,\theta) ~
           P_n^m (\cos\varphi) ~ \cos (m \theta)     \ ; \\
  b_{nm} = \frac {2n+1} {2 \pi (n+1)}
           \frac {(n-m) \factorial} {(n+m) \factorial} 
           \int_{r=1} d \Gamma ~ 
           h (\varphi,\theta) ~
           P_n^m (\cos\varphi) ~ \sin (m \theta)     \ ;
\label{eq-add200}
\end{eqnarray}
Once obtained the coefficients of the expansion, the additional
mass is obtained from
\begin{eqnarray} 
  A_{jj} && = \int_{r=1} d \Gamma ~ 
              \psi ~ \psi_{,r} \ ;  \\
         && = \int_{\varphi=0}^{\pi} d \varphi 
              \left  \{ 
          a_{n0}^2 [P_n (\cos \varphi)]^2 + \sum_{m=1}^n 
         (a_{nm}^2 + b_{nm}^2) [P_n^m (\cos \varphi)]^2 
              \right \} \ ; 
\label{eq-add210}
\end{eqnarray}
where the orthogonality property of the spherical harmonics 
was taken into account, and
\begin{eqnarray} 
   A_{jj} && = \int_{r=1} d \Gamma ~ \psi ~ \psi_{,r} \\
          && = \int_{\varphi=0}^{\pi} d \varphi
               \left  \{ 
               a_{n0}^2 [P_n (\cos \varphi)]^2 + 
               \sum_{m=1}^n (a_{nm}^2 + b_{nm}^2) 
               [P_n^m (\cos \varphi) ]^2 
               \right \} \ ;
\label{eq-add220}
\end{eqnarray}
and using the properties of the Legendre polynomials
\begin{eqnarray} 
   A_{jj} = \sum_{n=1}^{\infty} 
            \frac 2 {2n+1} \left  [ 
               a_{n0}^2 + \sum_{m=1}^n 
              (a_{nm}^2 + b_{nm}^2) 
              \frac {(n-m) \factorial}
                    {(n+m) \factorial} 
                             \right ] \ .
\label{eq-add230}
\end{eqnarray}
\section{Hemisphere in heave at very-low frequencies}
In this case the load $h(\varphi,\theta)=|\cos \varphi|$, so
\begin{equation} 
   a_n = {2n+1 \over 4 \pi(n+1)} 
         \int_{-1}^1 d \mu ~ | \mu | ~ P_n (\mu) \ ;
\label{eq-add240}
\end{equation}
as the $P_k$ are even (odd) for $k$ even (odd), only 
remains the even terms and then
\begin{equation} 
 a_n = {2n+1 \over 2\pi(n+1)} 
       \int_0^1 \mu ~ d \mu ~ P_n(\mu) \qquad \hbox {for $n$ even.}
\label{eq-add250}
\end{equation}
For computing the integral, the $P_n$ terms are generating 
in a recursive way from $P_0=1$, $P_1=\mu$ and $P_2$, ..., $P_n$ 
are obtained solving
\begin{equation} 
   ( n + 1) ~       P_{n+1} - 
   (2n + 1) ~ \mu ~ P_{n  } (\mu) + 
          n ~     ~ P_{n-1} (\mu) = 0 \ . 
\label{eq-add270}
\end{equation}
The coefficients of the polynomials $\mu P_n (\mu)$ are obtained
from the $P_n$ ones, and the integral is made in a direct way.
The final result is 
\begin{equation} 
   a_{33} = 1.7403 ~ \rho ~ R^3 \ ;
\label{eq-add280}
\end{equation}
corresponding to $a_{33}' = 0.83093$, that is, the non-dimensional 
coefficient respects to the hemisphere mass $2/3\pi\rho R^3$. 
\section{Hemisphere in heave at very-high frequencies}
In this case, the factor $\cos \theta$ made that the only 
non-null coefficients are the $a_{n1}$ ones. For obtaining
these, an integral from $\mu=1$ to $-1$ must be done, with
a function that contents the $P_n^1$ ones. These terms 
having a factor $\sqrt{1-\mu^2}$, so it is convenient 
to compute the integral in a numerically way. Our final 
numerical result is
\begin{equation} 
   a_{11} (\omega \to \infty) = 0.29806 ... 
                        \qquad (0.14231 ...) \ ; 
\label{eq-add300}
\end{equation}
where the value between parentheses is referred, always, to the
non-dimensional value respect to displaced mass.

\noindent
In brief, our estimates for the asymptotic values at very-low 
and very-high frequencies (subindex 1 for the surge mode and 3 
for the heave one), of the non-dimensional coefficients respects 
to the hemisphere mass $2/3\pi\rho R^3$ (denoted with primes) are
\begin {equation}
   \cases { 
     A'_{11} \to 0.5       & for $ K a \to 0      $ ; \cr
     A'_{11} \to 0.14231   & for $ K a \to \infty $ ; \cr
     A'_{33} \to 0.83093   & for $ K a \to 0      $ ; \cr
     A'_{33} \to 0.5       & for $ K a \to \infty $ . 
          }
\label{eq-nobl200}
\end{equation}
On the other hand, some literature values found for the surge/sway 
mode, e.g. see Sierevogel \cite{rf:sierevogel}, Prins \cite{rf:prins} 
(where only the intervals [0.25, 1.50] and [0.6, 1.5] are considered, 
respectively, so the extrapolations are rather doubtful), and for the 
heave one, e.g. see Korsmeyer \cite{rf:korsmeyer1} and Liapis 
\cite{rf:liapis}, are
\begin {equation}
   \cases { 
     A'_{11} \to 0.50   & for $ K a \to 0      $ ; \cr
     A'_{11} \to 0.25   & for $ K a \to \infty $ ; \cr
     A'_{33} \to 0.80   & for $ K a \to 0      $ ; \cr
     A'_{33} \to 0.45   & for $ K a \to \infty $ ; 
          }
\label{eq-nobl220}
\end{equation}
and the Hulme's \cite{rf:hulme} ones are
\begin{equation}
   \cases { 
 A'_{11}=0.5        & for $Ka \to 0     $ ; \cr
 A'_{11}=0.273~239..& for $Ka \to \infty$ ; \cr
 A'_{33}=0.830~951..& for $Ka \to 0     $ ; \cr
 A'_{33}=0.5        & for $Ka \to \infty$ ; 
          }
\label{eq-nobl250}
\end{equation}
Korsmeyer used a panel method with Fourier transform and 
complex impedance extended to very slow frequencies, while the 
Hulme's numerical results are obtained by spherical harmonics.
The Sierevogel, Prins and Liapis results are obtained with a 
panel method and Kelvin source. In general, the concordance 
between our estimates and the literature ones is good, except 
in our surge mode coefficient $A'_{11}$, which is about a half.

\section{Conclusions} 

We have considered a semi-numerical computation by spherical harmonics
for added mass coefficients of an oscillating hemi-sphere at very low 
and very high frequencies. The method is based by a physics insight 
and the computation is near-exact, so the results can be used to test 
the accuracy achieved by numerical codes, as finite-element and 
boundary-element ones.

\smallskip
\ltitle {Acknowledgment}
\smallskip

\noindent
This work was partially performed with the Free Software 
Foundation/GNU-Project resources, as GNU/Linux~OS and GNU/Octave,
and supported through grants CONICET-PIP-198/98 ({\it Germen-CFD}), 
SECyT-FONCyT-PICT-51 ({\it Germen}) and 
CAI+D-UNL-94-004-024.

\bibliography{navales}

\end{document}